\documentclass[a4paper,11pt]{article}
\usepackage{amssymb}
\usepackage{mathrsfs}
\usepackage{indentfirst}
\usepackage{amsmath}
\usepackage{amsfonts,amssymb}
\usepackage{bm}\usepackage{amsmath}
\usepackage{amssymb}
\usepackage{amsthm}
\usepackage{graphicx,color}
\catcode`\@=11 %\@addtoreset{equation}{section}

\catcode`\@=12

\newtheorem{Definition}{Definition}[section]
\newtheorem{Lemma}{Lemma}[section]

\newtheorem{Remark}{Remark}[section]

\title{On light-like extremal surfaces in curved
spacetimes\footnote{This work was supported by the NNSF of China
(Grant Nos. 11026151, 11101001) and the Anhui Provincial
University's Natural Science Foundation (Grand No. KJ2010A130).}}

\author{Shou-Jun Huang\;\;and\;\;Chun-Lei He\footnote{Corresponding author. Tel: +86 553 4810892; E-mail address: hcl026@126.com (Chun-Lei He)}\\
Department of Mathematics, Anhui Normal University\\ Wuhu 241000,
China }\date{ }
\begin{document}
\maketitle
\begin{abstract}In this paper, we are concerned with light-like extremal
surfaces in curved spacetimes.  It is interesting to find that under
a diffeomorphic transformation of variables, the light-like extremal
surfaces can be described by a system of nonlinear geodesic
equations. Particularly, we investigate the light-like extremal
surfaces in Schwarzschild spacetime in detail and some new special
solutions are derived systematically with aim to comparing with the
known results and illustrating the method. \vspace{3mm}


 \noindent{\bf Key words and phrases}: light-like extremal surfaces, geodesic equations,
  curved spacetimes, Schwarzschild spacetime.
%\vspace{3mm}

%\noindent{\bf 2000 Mathematics Subject Classification}:
\end{abstract}
\newpage
\section{Introduction}

In recent years the string and membrane theory has drawn great
interest. The reason lies in that this theory is not only a possible
unification model, but also it tightly relates to extremal surfaces
in physical spacetimes from mathematical point of view. For the
theory of extremal surfaces in flat Minkowski space-time, there have
been many deep results such as Calabi \cite{c} and Cheng and Yau
\cite{cy} for space-like case, the papers \cite{ac}, \cite{ba},
\cite{ksz}, \cite{m} for time-like case, and \cite{g3}-\cite{g4} for
the case of mixed type.  Hoppe {\it{et al.}} \cite{bh} derive the
equation for a classical relativistic membrane moving in Minkowski
space $\mathbb{R}^{1+3}$ and give some classical solutions, see also
\cite{h}. Lindblad \cite{lin} studies this equation and proves the
global existence of smooth solution for small initial data.

Recently, in the papers \cite{hk}-\cite{kzz} Kong {\it{et al.}}
restudy the dynamics of relativistic string in Minkowski space
$\mathbb{R}^{1+n}$. Based on the geometric properties enjoyed by the
extremal surface equations, they give a sufficient and necessary
condition for global existence of extremal surfaces without
space-like point in $\mathbb{R}^{1+n}$ with given initial data.
Moreover, numerical analysis shows that various topological
singularities will develop in finite time in the motion of a string.
Surprisingly, they obtain a general solution formula for the highly
coupled nonlinear equations. With the aid of this solution formula,
they succeed in  proving that the motion of a closed string is
always time-periodic. For the theory of extremal surfaces in curved
spacetimes, there are very few results so far (see e.g., \cite{g1}).
For a relativistic string, He and Kong \cite{hek1}-\cite{hek}
analyze the basic governing equations (see (\ref{1}) or (\ref{2})
below) and their inherent interesting properties. Under suitable
assumptions, they also provide some positive results on the global
existence for the motion of relativistic strings in  a general
curved spacetime.

 For  light-like extremal surfaces, it is well-known that the induced
metrics are degenerate by definition such that it is not clear how
to define the mean curvature and what kind of light-like surfaces
can be treated as extremal (cf. \cite{go}, \cite{he}). However, we
would like to point out that the system (\ref{2}) below is still
valid to describe  light-like surfaces. Kong {\it et al.} \cite{kzz}
give some remarks on light-like extremal surfaces in flat Minkowski
space. In \cite{hk}, Huang and Kong  consider three-dimensional
light-like extremal sub-manifolds in Minkowski spacetime and present
some explicit examples. Gorkavyy \cite{go} distinguishes two
particular classes of light-like surfaces in Minkowski spacetime by
applying one particular deformability property of surfaces. He
\cite{he} studies the light-like extremal surfaces in flat Minkowski
spacetime $\mathbb{R}^{1+(1+n)}$ and gives a necessary and
sufficient condition to obtain  explicit solution formulas. In
physics, another topic related to light-like surfaces is the
so-called {\it null string}, which was  introduced by Schild
\cite{schild} (see also \cite{kl}) and later developed by Dabrowski
and Larsen \cite{dl}  in curved spacetimes. In \cite{dl}, the
equations for a null string and the associated constraints in
Schwarzschild spacetime are derived and the dynamics for some
special solutions are also discussed.

This paper is devoted to investigating  light-like extremal surfaces
in a general curved spacetime. Here by light-like extremal surfaces,
we mean that for such surfaces, they not only satisfy the
light-likeness assumption (see (\ref{6}) below), but also they can
be described by the equations (\ref{2}). It is the light-likeness
condition that we are able to simplify the governing equations
(\ref{2}) into a system of geodesic equations in curved spacetimes.
Then based on the geodesic equations, we particularly study the
light-like extremal surfaces  in Schwarzschild spacetime and derive
some special solutions to illustrate the method in the present
paper.

The reminder of paper is organized as follows. In Section 2, we
introduce the equations for light-like extremal surfaces in curved
spacetimes. By the light-likeness assumption and a diffeomorphic
transformation of variables, we succeed in simplifying these
equations and obtaining  a system of geodesic equations. Section 3
is devoted to concerning the light-like extremal surfaces in
Schwarzschild spacetime. Many examples are provided systematically,
which have verified the method presented in this paper for studying
the light-like extremal surfaces. Finally, conclusions and some
important remarks are given in Section 4.

\section{Light-like extremal surfaces in curved spacetimes}
In this section, we concern the equations for light-like extremal
surfaces in a general curved spacetime $(\mathscr{N}, \tilde{g})$,
which is a  Lorentzian manifold. Some properties enjoyed by these
equations are also discussed.

For a two-dimensional  extremal surface, denoted by $S$, the local
coordinates are supposed to be $(\zeta^0,\zeta^1)$ and sometimes for
simplicity, we denote $\zeta^0=t,\; \zeta^1=\theta$. If the
 extremal surface is time-like, then the corresponding Euler-Lagrange
equations read
\begin{equation}\label{1}g^{\alpha\beta}\left(x^{\mu}_{\alpha\beta}+\tilde{\Gamma}^{\mu}_{
\lambda\nu}x^{\lambda}_{\alpha}x^{\nu}_{\beta}\right)=0,\quad(\mu=0,1,\cdots,n)\end{equation}
where $x^{\mu}_{\alpha}=\frac{\partial x^{\mu}}{\partial
\zeta^{\alpha}},x^{\mu}_{\alpha\beta}=\frac{\partial^2
x^{\mu}}{\partial \zeta^{\alpha}\partial\zeta^{\beta}}$,
$\alpha,\beta=0, 1,~\lambda, \mu, \nu, \rho=0,1,\cdots, n$.
$g^{\alpha\beta}$ is the inverse of the induced metric
$g_{\alpha\beta}=\tilde{g}_{\mu\nu}x^{\mu}_{\alpha}x^{\nu}_{\beta}$
on the extremal surface $S$. The coordinates
$x(t,\theta)=(x^0(t,\theta),x^1(t,\theta),\cdots, x^n(t,\theta))^T$
describe the surface $S$ and $\tilde{\Gamma}^{\mu}_{\lambda\nu}$
stand for the connections of the ambient metric $\tilde{g}$. The
above equations  can be rewritten in the following form
\begin{equation}\label{2}g_{11}x^{\mu}_{tt}-2g_{01}x^{\mu}_{t\theta}+g_{00}x^{\mu}_{\theta\theta}+g_{11}
\tilde{\Gamma}^{\mu}_{\nu\rho}x_t^{\nu}x_t^{\rho}
-2g_{01}\tilde{\Gamma}^{\mu}_{\nu\rho}x_t^{\nu}x_{\theta}^{\rho}+g_{00}\tilde{\Gamma}^{\mu}_{\nu\rho}x_{\theta}^{\nu}x_{\theta}^{\rho}=0,
\end{equation}
see for example \cite{hek} for details. Based on geometric
properties for (\ref{2}), He and Kong \cite{hek} prove a small-data
global result for a  string moving in general curved spacetimes.

 In this paper, we consider the Cauchy problem for (\ref{2}) with the following initial data
\begin{equation}\label{3}t=0:\quad x^{\mu}=\varphi^{\mu}(\theta),\quad
x^{\mu}_{t}=\psi^{\mu}(\theta)\;\quad(\mu=0,1,\cdots,n),
\end{equation}
where $\varphi^{\mu}(\theta)$ are $C^2$-smooth functions with
bounded $C^2$ norm, while $\psi^{\mu}(\theta)$ are $C^1$-smooth
functions with bounded $C^1$ norm. In physics,
$\varphi(\theta)=(\varphi^0(\theta),\varphi^1(\theta),\cdots,\varphi^n(\theta))$
and
$\psi(\theta)=(\psi^0(\theta),\psi^1(\theta),\cdots,\psi^n(\theta))$
stand for the initial position and initial velocity of a
relativistic string, respectively.  Moreover, $\varphi(\theta)$ and
$\psi(\theta)$ satisfy the following light-likeness condition
\begin{equation}\label{4}(g_{01}[\varphi,\psi](\theta))^2-g_{00}[\varphi,\psi](\theta)g_{11}[\varphi,\psi](\theta)\equiv0,
\end{equation}
in which
$$g_{00}[\varphi,\psi](\theta)\overset{\text{def}}{=}\tilde{g}_{\mu\nu}(\varphi)\psi^{\mu}\psi^{\nu},\;\;
g_{01}[\varphi,\psi](\theta)\overset{\text{def}}{=}\tilde{g}_{\mu\nu}(\varphi)\psi^{\mu}\varphi^{\nu}_{\theta},\;\;
g_{11}[\varphi,\psi](\theta)\overset{\text{def}}{=}\tilde{g}_{\mu\nu}(\varphi)\varphi^{\mu}_{\theta}\varphi^{\nu}_{\theta}.$$

Now we denote
\begin{equation}\label{5}\Delta=\Delta(t,\theta)\overset{\text{def}}{=} g_{01}^2-g_{00}g_{11},
\end{equation}
and introduce the following definition for light-like extremal
surfaces, which should be due to Kong ${\it et\; al.}$ \cite{kzz}
essentially.
\begin{Definition}Given an ambient curved spacetime $(\mathscr{N}, \tilde{g})$, a surface $S$ described by
 $x(t,\theta)=(x^0(t,\theta),x^1(t,\theta),\cdots,x^n(t,\theta))^T$ is called
 light-like extremal,
 if $x(t,\theta)$ satisfies the equations (\ref{2}) and the following
 light-likeness
 assumption
 \begin{equation}\label{6}\Delta\equiv0.
\end{equation}
\end{Definition}
\begin{Remark}In fact, for a surface $S$ in curved spacetimes, if $\Delta>0$ at every point in $S$, then
the surface is said to be entire time-like; if $\Delta<0$ at every
point in $S$, then the surface is entire space-like; if
$\Delta\equiv0$ at every point in $S$, then the surface is said to
be entire light-like; if a connected surface contains both a
time-like part and a space-like part simultaneously, then it is
 of mixed type.
\end{Remark}
\noindent By the above remark, we can see that the meaning of
light-like extremal surfaces here requires the equations (\ref{2})
should hold additionally on the light-like surfaces. This method to
handle light-like surfaces shall  bring us much more convenience
 in the following discussion.

\begin{Remark}By making use of the
system (\ref{2}) and the formulas for
$\tilde{\Gamma}^{\mu}_{\nu\rho}$, straightforward computations show
that\begin{equation}\label{8}g_{11}\Delta_t-g_{01}\Delta_{\theta}=2\Delta\left(\frac{\partial
g_{11}}{\partial t}-\frac{\partial g_{01}}{\partial \theta}\right),
\end{equation}
which implies the compatibility between the system (\ref{2}) and the
light-likeness condition (\ref{6}).\end{Remark}

Introduce
\begin{equation}\label{9}\lambda(t,\theta)\overset{\text{def}}{=}-\frac{g_{01}}{g_{11}},
\end{equation}
then the system (\ref{2}) can be rewritten in the following form
\begin{equation}\label{16}x^{\mu}_{tt}+2\lambda
x^{\mu}_{t\theta}+\lambda^2x_{\theta\theta}^{\mu}+\tilde{\Gamma}_{\nu\rho}^{\mu}(x)
\left(x_t^{\nu}x_t^{\rho}+2\lambda
x_t^{\nu}x_{\theta}^{\rho}+\lambda^2x_{\theta}^{\nu}x_{\theta}^{\rho}\right)=0.\end{equation}
First, we have
\begin{Lemma}Under the assumption (\ref{6}), we claim that $\lambda=\lambda(t,\theta)$
satisfies the following Burgers equation
\begin{equation}\label{10}\lambda_t+\lambda \lambda_{\theta}=0
\end{equation}
on the existence domain of smooth solution $x(t,\theta)$.
\end{Lemma}
{\bf Proof.} By the light-likeness condition (\ref{6}), direct
computation gives
\begin{eqnarray}\label{11}&&\lambda_t+\lambda \lambda_{\theta}\nonumber\\
&=&\lambda_t+\left(\frac{\lambda^2}{2}\right)_{\theta}=\left(-\frac{g_{01}}{g_{11}}\right)_t+\left(\frac12\frac{g_{00}}{g_{11}}\right)_{\theta}\nonumber\\
&=&-\frac{1}{g_{11}^2}\left(\frac{\partial g_{01}}{\partial
t}g_{11}-g_{01}\frac{\partial g_{11}}{\partial
t}\right)+\frac{1}{2g_{11}^2}\left(g_{11}\frac{\partial
g_{00}}{\partial\theta}-g_{00}\frac{\partial
g_{11}}{\partial\theta}\right)\nonumber\\
&=&\frac{1}{g_{11}^2}\left[g_{01}\left(\frac{\partial\tilde{g}_{\mu\nu}}{\partial
x^{\rho}}x_t^{\rho}x_{\theta}^{\mu}x_{\theta}^{\nu}+2\tilde{g}_{\mu\nu}x_{t\theta}^{\mu}x_{\theta}^{\nu}\right)\right.\nonumber\\
&&\left. -g_{11}\left(\frac{\partial\tilde{g}_{\mu\nu}}{\partial
x^{\rho}}x_t^{\rho}x_{t}^{\mu}x_{\theta}^{\nu}+\tilde{g}_{\mu\nu}x_{tt}^{\mu}x_{\theta}^{\nu}+\tilde{g}_{\mu\nu}x_{t}^{\mu}x_{\theta
t}^{\nu}\right)
\right.\nonumber\\
&&\left.+\frac12g_{11}\left(\frac{\partial\tilde{g}_{\mu\nu}}{\partial
x^{\rho}}x_{\theta}^{\rho}x_{t}^{\mu}x_{t}^{\nu}+2\tilde{g}_{\mu\nu}x_{t\theta}^{\mu}x_{t}^{\nu}\right)\right.\nonumber\\
&&\left.
-\frac12g_{00}\left(\frac{\partial\tilde{g}_{\mu\nu}}{\partial
x^{\rho}}x_{\theta}^{\rho}x_{\theta}^{\mu}x_{\theta}^{\nu}+2\tilde{g}_{\mu\nu}x_{\theta\theta}^{\mu}x_{\theta}^{\nu}\right)\right]\nonumber\\
&=&-\frac{1}{g_{11}}\left[\tilde{g}_{\mu\nu}x_{\theta}^{\nu}\left(x^{\mu}_{tt}+2\lambda
x^{\mu}_{t\theta}+\lambda^2x^{\mu}_{\theta\theta}\right)+\frac{\partial\tilde{g}_{\mu\nu}}{\partial
x^{\rho}}x_t^{\mu}\left(x_t^{\rho}x^{\nu}_{\theta}-\frac12x^{\rho}_{\theta}x^{\nu}_{t}\right)
\right.\nonumber\\
&&\left.+\lambda\frac{\partial\tilde{g}_{\mu\nu}}{\partial
x^{\rho}}x_t^{\rho}x^{\mu}_{\theta}x_{\theta}^{\nu}+\frac{\lambda^2}{2}\frac{\partial\tilde{g}_{\mu\nu}}{\partial
x^{\rho}}x_{\theta}^{\rho}x^{\mu}_{\theta}x_{\theta}^{\nu}\right].
\end{eqnarray}
By utilizing (\ref{2})  it follows from (\ref{11}) that
\begin{eqnarray*}&&\lambda_t+\lambda \lambda_{\theta}\nonumber\\
&=&-\frac{1}{g_{11}}\left[-x_{\theta}^{\nu}\tilde{g}_{\mu\nu}\left(\tilde{\Gamma}_{\lambda\rho}^{\mu}x_t^{\lambda}x_t^{\rho}+2\lambda\tilde{\Gamma}^{\mu}
_{\lambda\rho}x_t^{\lambda}x_{\theta}^{\rho}+\lambda^2\tilde{\Gamma}_{\lambda\rho}^{\mu}x_{\theta}^{\lambda}x_{\theta}^{\rho}\right)\right.\nonumber\\
&&\left.+\frac{\tilde{g}_{\mu\nu}}{\partial
x^{\rho}}x_t^{\mu}\left(x_t^{\rho}x_{\theta}^{\nu}-\frac12x_{\theta}^{\rho}x_t^{\nu}\right)+\lambda
\frac{\partial\tilde{g}_{\mu\nu}}{\partial
x^{\rho}}x_t^{\rho}x_{\theta}^{\mu}x_{\theta}^{\nu}+\frac{\lambda^2}{2}\frac{\partial\tilde{g}_{\mu\nu}}{\partial
x^{\rho}}x_{\theta}^{\rho}x_{\theta}^{\mu}x_{\theta}^{\nu}\right]\nonumber\\
&=&-\frac{1}{g_{11}}\left[-\frac12x_{\theta}^{\nu}\left(\frac{\partial\tilde{g}_{\nu\rho}}{\partial
x^{\lambda}}+\frac{\partial\tilde{g}_{\nu\lambda}}{\partial
x^{\rho}}-\frac{\partial\tilde{g}_{\lambda\rho}}{\partial x^{\nu}}
\right)\left(x_t^{\lambda}x_t^{\rho}+2\lambda x_t^{\lambda}x_{\theta}^{\rho}+\lambda^2x_{\theta}^{\lambda}x_{\theta}^{\rho}\right)\right.\nonumber\\
&&\left.+\frac{\tilde{g}_{\mu\nu}}{\partial
x^{\rho}}x_t^{\mu}\left(x_t^{\rho}x_{\theta}^{\nu}-\frac12x_{\theta}^{\rho}x_t^{\nu}\right)+\lambda
\frac{\partial\tilde{g}_{\mu\nu}}{\partial
x^{\rho}}x_t^{\rho}x_{\theta}^{\mu}x_{\theta}^{\nu}+\frac{\lambda^2}{2}\frac{\partial\tilde{g}_{\mu\nu}}{\partial
x^{\rho}}x_{\theta}^{\rho}x_{\theta}^{\mu}x_{\theta}^{\nu}\right]\nonumber\\
&=&0.
\end{eqnarray*}
Thus, the proof is completed.   $\quad\quad\blacksquare$

\begin{Remark}We would like to point out that the method in \cite{hek} can not be directly applied here due
to the light-likeness condition $\Delta=0$ (see Theorem 2.1 in
\cite{hek}).
\end{Remark}
We consider the Burgers equation (\ref{10}) associated with the
following initial data
\begin{equation}\label{12}t=0:\quad
\lambda=\Lambda(\theta)\overset{\text{def}}{=}-\frac{g_{01}[\varphi,\psi](\theta)}{g_{11}[\varphi,\psi](\theta)}.
\end{equation}
It is well-known that in order the global existence of smooth
solutions for the Cauchy problem (\ref{10}) and (\ref{12}) exists,
we can impose the following sufficient and necessary condition on
the initial data (cf. \cite{li}):
\begin{equation}\label{13}\Lambda'(\theta)\geqslant0, \quad
\forall\; \theta\in\mathbb{R}.
\end{equation}
Under the assumption (\ref{13}), the smooth solution
$\lambda(t,\theta)$ then can be solved as
\begin{equation}\label{14}\lambda(t,\theta)=\Lambda(\vartheta(t,\theta)),\end{equation}
where $\vartheta(t,\theta)$ is the inverse function of
\begin{equation}\label{15}\theta=\vartheta+\Lambda(\vartheta)t
\end{equation}
for any fixed $ t\geqslant0$.

Now we are ready to introduce the following transformation of the
variables
\begin{equation}\label{17}(t,\theta)\rightarrow
(t,\vartheta),
\end{equation}
where $\vartheta=\vartheta(t,\theta)$ is defined through the
equation (\ref{15}) implicitly.

 We have
\begin{Lemma}Under the assumption (\ref{13}), the mapping defined by (\ref{17}) is globally diffeomorphic.\end{Lemma}
{\bf Proof.}
 It is obvious that the mapping (\ref{17}) defined through (\ref{15}) is
well-defined on $\mathbb{R}^+\times\mathbb{R}$. Moreover, under the
assumption (\ref{13}), we have
\begin{eqnarray}\label{19}\mathscr{J}&\overset{\text{def}}{=}&\frac{\partial(t,\vartheta)}{\partial(t,\theta)}=\left|\begin{array}{cc}1&0\\
\frac{\partial\vartheta}{\partial t}&
\frac{\partial\vartheta}{\partial
\theta}\end{array}\right|=\frac{\partial\vartheta(t,\theta)}{\partial
\theta}=\frac{1}{1+\Lambda'(\vartheta)t}>0,
\end{eqnarray} for every
$(t,\theta)\in\mathbb{R}^{+}\times\mathbb{R}$.  Then by the
Hadamard's Lemma \cite{gordon}, we can conclude that the mapping
(\ref{17}) is globally diffeomorphic. Thus, we complete the proof.
$\quad\quad\blacksquare$

 Furthermore, direct calculation shows that under
the new coordinates $(t,\vartheta)$, the system (\ref{16}) or
(\ref{2}) can be equivalently reduced into the following form,
\begin{equation}\label{21}y^{\mu}_{tt}+\tilde{\Gamma}_{\nu\rho}^{\mu}(y)y_t^{\nu}y_{t}^{\rho}=0,\quad
(\mu=0,1,\cdots,n)
\end{equation}
where $y^{\mu}=y^{\mu}(t,\vartheta)=x^{\mu}(t,\theta(t,\vartheta))$.
Due to $\vartheta=\theta$ at $t=0$, the initial data for the new
system (\ref{21}) take the following form
\begin{equation}\label{20}t=0:\quad y^{\mu}=\varphi^{\mu}(\vartheta),\quad
y^{\mu}_{t}=\psi^{\mu}(\vartheta)\;\quad(\mu=0,1,\cdots,n),
\end{equation}
and the assumption (\ref{13}) is equivalent to
\begin{equation}\label{222}\Lambda'(\vartheta)\geqslant0.\end{equation}
In addition, it is easy to see that the light-likeness condition
(\ref{6}) can be preserved for this kind of diffeomorphic
transformation since we have
\begin{equation}\label{23}\Delta(t,\theta)=\Delta(t,\vartheta)\vartheta_{\theta}^2,\end{equation}
where $\vartheta_{\theta}>0$ under the assumption (\ref{13}) or
(\ref{222}).

 It is interesting  to see that
in the new coordinates $(t,\vartheta)$, the system (\ref{21})
 is independent of the variable $\vartheta$. This means that for any point $\vartheta$ in the initial curve,
 the trajectory of this point must be geodesic in the background spacetime.  In other words, if the initial data
(\ref{3}) satisfy the light-likeness condition (\ref{4}) and the
assumption (\ref{13}), all $t$-curves of the light-like extremal
surfaces should  be geodesic. However, in general cases the smooth
solutions to (\ref{21}) can not exist globally in time due to the
appearance of nonlinearity arising  from the ambient curved
spacetimes. For the theory on geodesic equations in curved
spacetimes, we refer to the monograph by Chandrasekhar \cite{chan}.

\begin{Remark}If the background spacetime is flat,  then the
ambient connections vanish and the system (\ref{21}) will go back to
the equations studied in \cite{he}.
\end{Remark}
\section{Light-like extremal surfaces in Schwarzschild spacetime}
In this section, we  study the light-like extremal surfaces in
Schwarzschild spacetime, which is stationary, spherically symmetric
and asymptotically flat. In the spherical coordinates
$(\tau,r,\alpha,\beta)$, the Schwarzschild metric $\tilde{g}$ reads
\begin{equation}\label{22}ds^2=-\left(1-\dfrac{2m}{r}\right)d\tau^2+\left(1-\dfrac{2m}{r}\right)^{-1}dr^2+r^2
\left(d\alpha^2+\sin^2\alpha d\beta^2\right),\end{equation} where
$m$ is a positive constant standing for the universe mass.

For the above Schwarzschild metric, the equations of motion for
light-like extremal surfaces (\ref{21}) reduce to the following form
\begin{equation}\label{24}\left\{\begin{array}{l}\displaystyle{\tau_{tt}+\frac{2m}{r(r-2m)}\tau_t\,r_t=0,}\vspace{2mm}\\
\displaystyle{r_{tt}-(r-2m)\sin^2\alpha\,\beta_t^2-(r-2m)\,\alpha^2_t-\frac{m}{r(r-2m)}\,r_t^2+\frac{m(r-2m)}{r^3}\,\tau_t^2=0,}\vspace{2mm}\\
\displaystyle{\alpha_{tt}+\frac2r\,r_t\,\alpha_t-\frac12\sin2\alpha\,\beta^2_t=0,}\vspace{2mm}\\
\displaystyle{\beta_{tt}+\frac2r\,r_t\,\beta_t+\frac{2\cos\alpha}{\sin\alpha}\,\alpha_t\,\beta_t=0.}\end{array}\right.
\end{equation}
The first and last equation in (\ref{24}) can be easily integrated,
i.e.,
\begin{equation}\label{25}\tau_t=\frac{E(\vartheta)}{1-\frac{2m}{r}},
\end{equation}
\begin{equation}\label{26}\beta_t=\frac{L(\vartheta)}{r^2\sin^2\alpha},
\end{equation}
where $E(\vartheta)$ and $L(\vartheta)$ are two smooth functions and
will be determined by the initial data. Substituting (\ref{26}) into
the third equation in (\ref{24}) yields
\begin{equation}\label{27}r^4\sin^2\alpha\;\alpha_t^2=-L^2(\vartheta)\cos^2\alpha+K(\vartheta)\sin^2\alpha,
\end{equation}
where $K(\vartheta)$ is a nonnegative function. Combining
(\ref{25})-(\ref{27}) and the second equation in (\ref{24}) leads to
the following equation for $r$:
\begin{equation}\label{28}r_{tt}-\frac{m}{r(r-2m)}r_t^2+\frac{mE^2(\vartheta)}{r(r-2m)}-\frac{r-2m}{r^4}\left(K(\vartheta)+L^2
(\vartheta)\right)=0.
\end{equation}
 The equations
(\ref{25})-(\ref{27}) are analogous in some manner to the equations
for the motion of null strings (see e.g. \cite{dl}). However,  it is
worth pointing out that we need $\Delta(\vartheta)=0$, instead of
$g_{00}=g_{01}=0$ in the null string theory, see also the following
discussions.

By the initial data (\ref{3}), the functions $E(\vartheta),
L(\vartheta)$ and $K(\vartheta)$ are in fact given  by
\begin{equation}\label{29}E(\vartheta)=\psi_0(\vartheta)\left(1-\frac{2m}{\varphi_1(\vartheta)}\right),\quad
L(\vartheta)=\psi_3(\vartheta)\varphi^2_1(\vartheta)\sin^2(\varphi_2(\vartheta))\end{equation}
and
\begin{equation}\label{30}
K(\vartheta)=\varphi_1^4(\vartheta)\left[\psi_2^2(\vartheta)+\psi^2_3(\vartheta)\sin^2(\varphi_2(\vartheta))\cos^2(\varphi_2(\vartheta))\right].
\end{equation}
 Here the lower
and upper indices are used interchangeably without ambiguity .

In order to solve the light-like extremal surfaces in Schwarzschild
spacetime, it suffices to consider the Cauchy problem
(\ref{25})-(\ref{30}) and (\ref{3}) under the assumption
(\ref{222}). By computation, we have
\begin{equation}\label{311}\Lambda(\vartheta)=
-\frac{-\left(1-\frac{2m}{\varphi_1}\right)\varphi_0'\psi_0+\left(1-\frac{2m}{\varphi_1}\right)^{-1}\varphi_1'\psi_1
+\varphi_1^2\varphi_2'\psi_2+\varphi_1^2\sin^2\varphi_2\;\varphi_3'\psi_3}{-\left(1-\frac{2m}{\varphi_1}\right)
\left(\varphi'_0\right)^2+\left(1-\frac{2m}{\varphi_1}\right)^{-1}\left(\varphi'_1\right)^2
+\varphi_1^2\left(\varphi'_2\right)^2+\varphi_1^2\sin^2\varphi_2\;\left(\varphi'_3\right)^2},\end{equation}
and
\begin{eqnarray}\label{312}\Delta(0,\theta)&=&\Delta(0,\vartheta)\nonumber\\
&=&\left(\psi_0\varphi_1'-\psi_1\varphi_0'\right)^2+\left(1-\frac{2m}{\varphi_1}\right)
\varphi_1^2\left(\psi_0\varphi_2'-\varphi_0'\psi_2\right)^2\nonumber\\
&&+\left(1-\frac{2m}{\varphi_1}\right)\varphi_1^2\sin^2\varphi_2\left(\psi_0\varphi_3'-\psi_3\varphi_0'\right)^2-
\left(1-\frac{2m}{\varphi_1}\right)^{-1}\varphi_1^2\left(\psi_1\varphi_2'-\psi_2\varphi_1'\right)^2\nonumber\\
&&-\left(1-\frac{2m}{\varphi_1}\right)^{-1}\varphi_1^2\sin^2\varphi_2\left(\psi_1\varphi_3'-\psi_3\varphi_1'\right)^2
-\varphi_1^4\sin^2\varphi_2\left(\psi_3\varphi_2'-\psi_2\varphi_3'\right)^2,
\end{eqnarray}
which will be useful in the following discussions.

 Of course, this kind of Cauchy problem looks
very complicated and can not be solved  generally. So in what
follows we would like to investigate some particular solutions to
illustrate our method.

{\bf Example 1}: We select the initial data as follows
\begin{equation}\label{32}t=0:\quad\varphi=(\tau_0,r_0,\alpha_0(\vartheta),\vartheta),\quad
\psi=\left(\pm\left(1-\frac{2m}{r_0}\right)^{-1}r_1,r_1,0,0\right),
\end{equation}
where $\tau_0, r_0, r_1$ are constants and $\alpha_0(\vartheta)$ is
an arbitrary function of $\vartheta$. Further we suppose that
$r_0>2m$. For this kind of initial data, direct computation shows
that the assumptions (\ref{222}) are satisfied and
$\Delta(0,\vartheta)=0$. By (\ref{29})-(\ref{30}) and (\ref{32}), it
is easy to see that
\begin{equation}\label{33}L(\vartheta)=0,\;\;K(\vartheta)=0,\;\;E(\vartheta)=\pm\,
r_1.
\end{equation}
 Then it follows from (\ref{26})-(\ref{27}) and
(\ref{32})-(\ref{33}) that
\begin{equation}\label{34}\alpha=\alpha_0(\vartheta),\quad \beta=\vartheta.
\end{equation}
Moreover, we can obtain from (\ref{28}) the following equation for
$r$:
\begin{equation}\label{35}r_{tt}-\frac{m}{r(r-2m)}r_t^2+\frac{m\,r_1^2}{r(r-2m)}=0,
\end{equation} where we have made use of (\ref{32})
and (\ref{33}). It can be observed  that the solution to (\ref{35})
is given by
\begin{equation}\label{36}r=\pm\, r_1\,t+r_0.
\end{equation}
Thus, by integration we derive from (\ref{25}) and (\ref{36}) the
following solution
\begin{equation}\label{37}
r-r_0+2m\ln\frac{r-2m}{r_0-2m}=\pm\,(\tau-\tau_0),\quad\alpha=\alpha_0(\vartheta),\quad
\beta=\vartheta.
\end{equation}
%Now we remark that when $\sin \alpha_0(\vartheta)=0$, i.e.,
%$\alpha_0(\vartheta)=0$ or $\pi$, the solution (\ref{37})  still can
%be a solution to the system (\ref{24}) with the initial data
%(\ref{32}).
 Here we would like to emphasize that
the solution (\ref{37}) is a generalization of the so-called {\it
cone strings} in the null string theory, since the cone strings
require the coordinate $\alpha$ is a constant (see e.g., \cite{dl}).
$\qquad\qquad \blacksquare$

The previous discussions imply that for the initial data (\ref{32})
satisfying the condition (\ref{222}), the light-like extremal
surface is uniquely given by (\ref{37}). Of course, by the method in
the present paper, one can choose other various initial data to
construct the corresponding light-like extremal surfaces.

Now we do some manipulations for our problem before we construct
other types of light-like extremal surfaces. For simplicity we
assume that the initial datum
\begin{equation}\label{38}\psi_3(\vartheta)\equiv0,\end{equation}
which implies
\begin{equation}\label{39}L(\vartheta)\equiv0.\end{equation}
By (\ref{26}), it leads to
\begin{equation}\label{40}\beta=\varphi_3(\vartheta)\end{equation}
and the equation (\ref{27}) becomes
\begin{equation}\label{41}\alpha_t=\pm\frac{\sqrt{K}}{r^2},\end{equation}
where we have supposed that $\sin\alpha\neq0$.

  Let
\begin{equation}\label{42}z=r_t.\end{equation}
When $r_t\neq0$, it follows from the third equation in (\ref{28})
that
\begin{equation}\label{43}\frac{dz}{dr}=\frac{m}{r(r-2m)}z+\left[\frac{r-2m}{r^4}K-\frac{mE^2}{r(r-2m)}\right]
\frac{1}{z},\end{equation} which in turn implies that
\begin{eqnarray}\label{44}\frac{d}{dr}\left(\frac{z^2}{\frac{r-2m}{r}}\right)
=\frac{d}{dr}\left(-\frac{K}{r^2}+\frac{2mE^2}{r-2m}\right).\end{eqnarray}
Integrating (\ref{44}) and noting the initial data, we have
\begin{equation}\label{45}r_t^2=\left(C-\frac{K}{r^2}\right)\left(1-\frac{2m}{r}\right)+\frac{2m}{r}E^2,\end{equation}
where
\begin{equation}\label{46}C=\frac{\psi_1^2\varphi_1^3+K(\varphi_1-2m)-2mE^2\varphi_1^2}{\varphi_1^2(\varphi_1-2m)}.\end{equation}
 By considering $r$ as a function of $\alpha$ (instead of $t$), we
obtain from the equation  (\ref{41}) and (\ref{45}) the following
equation
\begin{equation}\label{47}\left(\frac{du}{d\alpha}\right)^2=2mu^3-u^2+2mAu+B\triangleq
g(u),\end{equation} where $u=\frac1r$ and
\begin{equation}\label{48}A=-\frac{1}{\varphi_1^2}+\frac{(\varphi_1-2m)^2\psi_0^2-\varphi_1^2\psi^2_1}{(\varphi_1-2m)\varphi_1^5\psi_2^2},\end{equation}
\begin{equation}\label{49}B=\frac{1}{\varphi_1^2}+\frac{-2m(\varphi_1-2m)^2\psi_0^2+\varphi_1^3\psi^2_1}{(\varphi_1-2m)\varphi_1^6\psi_2^2}.\end{equation}

 Once equation (\ref{47})
has been solved for $u(\alpha)$, the solution can be completed by
direct quadratures of the following equations
\begin{equation}\label{50}\frac{dt}{d\alpha}=\pm\frac{1}{\sqrt{K}u^2}\quad \text{and}\quad \frac{d\tau}{d\alpha}=\mp\frac{E}{\sqrt{K}u^2(2mu-1)}.\end{equation}

{\bf Example 2}:  We choose the following initial data,
\begin{equation}\label{51}t=0:\quad\varphi=(\tau_0,r_0,\alpha_0,\vartheta),\quad
\psi=\left(f(\vartheta),0,\pm\frac{1}{r_0^2}\sqrt{r_0(r_0-2m)f^2(\vartheta)},0\right),
\end{equation}
where $\tau_0, r_0, \alpha_0$ are constants and $ r_0>2m$,
$f(\vartheta)$ is a non-zero smooth function. For the above initial
data, it is easy to see that the condition (\ref{222}) and
$\Delta(0,\vartheta)=0$ can be satisfied. Meanwhile, by
(\ref{29})-(\ref{30}) and (\ref{51}), we can observe that
\begin{equation}\label{52}L(\vartheta)=0,\;\;K(\vartheta)=r_0(r_0-2m)f^2(\vartheta)>0,
\;\;E(\vartheta)=\left(1-\frac{2m}{r_0}\right)f(\vartheta).
\end{equation}
Moreover, by the initial data (\ref{51}) we observe that the
equation (\ref{47}) can be rewritten as
\begin{equation}\label{53}\left(\frac{du}{d\alpha}\right)^2=2mu^3-u^2+\frac{r_0-2m}{r_0^3}=g(u).\end{equation}
Obviously, the solution relates to the disposition of the roots of
the cubic equation $g(u)=0$. In fact, we have
\begin{equation}\label{54}g(u)=2m\left(u-u_1\right)(u-u_2)(u-u_3),\end{equation}
where
\begin{equation}\label{55}u_1=\frac{1}{r_0}, u_2=\frac{r_0-2m+\sqrt{(r_0-2m)(r_0+6m)}}{4mr_0},
u_3=\frac{r_0-2m-\sqrt{(r_0-2m)(r_0+6m)}}{4mr_0}.\end{equation} It
is noted from (\ref{55}) that the third root $u_3$ of $g(u)$ is
always negative.

Returning to equation (\ref{53}), we will study the solution in the
following several cases. Case I: $r_0=3m$, equivalently,
$u_1=u_2=\frac{1}{3m}$. The equation (\ref{53}) becomes
\begin{equation}\label{56}\left(\frac{du}{d\alpha}\right)^2=2m\left(u+\frac{1}{6m}\right)\left(u-\frac{1}{3m}\right)^2.\end{equation}
Obviously, $u=\frac{1}{r_0}=\frac{1}{3m}$ is a solution to
(\ref{56}) satisfying the initial data. Then by (\ref{50})  we can
easily obtain the following light-like extremal surface
\begin{equation}\label{57}\tau=f(\vartheta)t+\tau_0,\quad r=3m,\quad\alpha=\pm\frac{\sqrt{f^2(\vartheta)}}{3\sqrt{3}m}t+\alpha_0,\quad \beta=\vartheta.\end{equation}
This special solution has been discovered in  \cite{dl} and we will
not discuss it further on.

\noindent Case II: $2m<r_0<3m$, equivalently, $u_1>u_2$. Bearing in
mind that $u$ initially equals to $u_1$ and the domain of $u$ should
be chosen so that $g(u)\geqslant0$, we should consider the solution
in the range $u_1\leqslant u <\frac{1}{2m}$ in this case. We now
make the substitution
\begin{equation}\label{58}u=u_2+(u_1-u_2)\sec^2\frac{\xi}{2}.\end{equation}
It is noted from (\ref{58}) that $\xi=0$ when $u=u_1$. Substituting
(\ref{58}) into (\ref{53}) gives
\begin{equation}\label{59}\left(\frac{d\xi}{d\alpha}\right)^2=2m(u_1-u_3)\left(1-k^2\sin^2\frac{\xi}{2}\right),\quad
k^2=\frac{u_2-u_3}{u_1-u_3}\;\;(0<k^2<1).\end{equation}
Then $\alpha$ can be solved in terms of the Jacobian elliptic
integral,
\begin{equation}\label{61}\alpha=2\left[2m(u_1-u_3)\right]^{-\frac12}F\left(\frac{\xi}{2},k\right)+\alpha_0,\end{equation}
where
\begin{equation}\label{62}F(\chi,k)=\int_0^{\chi}\frac{d\gamma}{\sqrt{1-k^2\sin^2\gamma}}.\end{equation}
In the derivation of (\ref{61}), we have made use of the initial
data (\ref{51}).

 \noindent Case III: $r_0>3m$,
equivalently, $u_1<u_2$. Noting the initial data, we have to study
the solution in the range $0<u\leqslant u_1$ for the present case.
We can choose the following substitution
\begin{equation}\label{63}u=u_3+\frac12(u_1-u_3)(1-\cos\xi).\end{equation}
For this selection, we have $u=u_1$ when $\xi=\pi$. By this
substitution, the equation (\ref{53}) reduces to
\begin{equation}\label{64}\left(\frac{d\xi}{d\alpha}\right)^2=2m(u_2-u_3)\left(1-k^2\sin^2\frac{\xi}{2}\right),\quad
k^2=\frac{u_1-u_3}{u_2-u_3}\;\;(0<k^2<1).\end{equation}
Similarly, $\alpha$ can be expressed by the Jacobian elliptic
integral as
\begin{equation}\label{66}\alpha=2\left[2m(u_2-u_3)\right]^{-\frac12}F\left(\frac{\xi}{2},k\right)+C_0,\end{equation}
where $C_0$ is determined by the initial data in the following form
\begin{equation}\label{67}C_0=\alpha_0-2\left[2m(u_2-u_3)\right]^{-\frac12}F\left(\frac{\pi}{2},k\right).\end{equation}

Thus, we have completed the discussions for this example.
$\qquad\qquad\qquad\blacksquare$

{\bf Example 3}: We now turn to consider the following initial data
\begin{equation}\label{68}t=0:\quad\varphi=(\varphi_0(\vartheta),\varphi_1(\vartheta),\varphi_2(\vartheta),\beta_0),\quad
\psi=\left(\varphi_0'(\vartheta),\varphi_1'(\vartheta),\varphi_2'(\vartheta),0\right),
\end{equation}
which implies that
\begin{equation}\label{69}\Lambda(\vartheta)=-1,\quad \Delta(\vartheta)=0\end{equation}
and the condition (\ref{222}) are satisfied automatically.

By the definition of $\Lambda(\vartheta)$, we need to require
$g_{11}\neq0$, which yields that the coefficient of $u$ in
(\ref{47}) should be a non-zero constant. For brevity we assume that
the constant
\begin{equation}\label{70}B=\frac{1}{\varphi_1^2}+\frac{-2m(\varphi_1-2m)^2\psi_0^2+\varphi_1^3\psi^2_1}{(\varphi_1-2m)\varphi_1^6\psi_2^2}=0,\end{equation}
equivalently,
\begin{equation}\label{71}\left(\varphi_2'\right)^2=\frac{2m(\varphi_1-2m)^2(\varphi_0')^2-\varphi_1^3(\varphi_1')^2}{(\varphi_1-2m)\varphi_1^4},\end{equation}
where we have made use of the initial data (\ref{68}). Then $A$
becomes a non-zero constant and in turn the induced metric
components $g_{00}=g_{01}=g_{11}\neq0$.

\begin{Remark}The above argument shows that for some light-like extremal surfaces, each induced metric component
may be non-zero. This important feature is very different from that
of null strings, since for null strings, one requires
$g_{00}=g_{01}=0$ (see e.g., \cite{dl}).
\end{Remark}

Next we discuss the light-like extremal surface in detail. Here for
illustration purpose, we only consider the case that
\begin{equation}\label{72}\varphi_0=\vartheta,\quad \varphi_1=r_0\end{equation}
and then from (\ref{71}) we obtain $\varphi_2$ as
\begin{equation}\label{73}\varphi_2=\pm\frac{\sqrt{2m(r_0-2m)}}{r_0^2}\,\vartheta.\end{equation}
So the initial data now become
\begin{equation}\label{74}t=0:\quad\varphi=(\vartheta,r_0,\pm\frac{\sqrt{2m(r_0-2m)}}{r_0^2}\,\vartheta,\beta_0),\quad
\psi=(1,0,\pm\frac{\sqrt{2m(r_0-2m)}}{r_0^2},0).
\end{equation}
Meanwhile, we have
\begin{equation}\label{75}L=0,\quad E=1-\frac{2m}{r_0},\quad K=2m(r_0-2m)>0\end{equation}
and the equation (47) can be rewritten as
\begin{equation}\label{76}\left(\frac{du}{d\alpha}\right)^2=2mu^3-u^2+\frac{r_0-2m}{r_0^2}u=2mu\left(u-\frac{1}{r_0}\right)\left(u-\frac{r_0-2m}{2mr_0}\right).\end{equation}
%\begin{Remark}In the present situation, the induced metric components
%\begin{equation}\label{77}g_{00}=g_{01}=g_{11}=-\frac{(r_0-2m)^2}{r_0^2}\neq0.\end{equation}\end{Remark}

In a similar manner to Example 2, we divide the  arguments into
several cases according to the values of $r_0$.

 \noindent Case I: $r_0=4m$,
equivalently, $\frac{1}{r_0}=\frac{r_0-2m}{2mr_0}$. For this case,
we can easily obtain the following simple solution
\begin{equation}\label{78}\tau=t+\vartheta,\quad r=4m,\quad \alpha=\pm\frac{1}{8m}(t+\vartheta),\quad \beta=\beta_0.\end{equation}

 \noindent Case II: $r_0>4m$,
equivalently, $\frac{1}{r_0}<\frac{r_0-2m}{2mr_0}$. We only need to
investigate the solution in the range $0<u\leqslant\frac{1}{r_0}$.
Take the following substitution
\begin{equation}\label{79}u=\frac{1}{2r_0}(1-\cos\xi).\end{equation}
We note that $u=\frac{1}{r_0}$ when $\xi=\pi$. So the equation
(\ref{76}) reduces into
\begin{equation}\label{80}\left(\frac{d\xi}{d\alpha}\right)^2=\frac{r_0-2m}{r_0}\left(1-k^2\sin^2\frac{\xi}{2}\right),\quad
k^2=\frac{2m}{r_0-2m}\;\;(0<k^2<1)\end{equation} and then $\alpha$
can be solved by
\begin{equation}\label{81}\alpha=2\left(\frac{r_0}{r_0-2m}\right)^{\frac12}F\left(\frac{\xi}{2},k\right)+C_0,\end{equation}
where $C_0$ is the integration constant,
$$C_0=-2\left(\frac{r_0}{r_0-2m}\right)^{\frac12}F\left(\frac{\pi}{2},k\right)\pm\frac{\sqrt{2m(r_0-2m)}}{r_0^2}\vartheta.$$

 \noindent Case III: $2m<r_0<4m$,
equivalently, $\frac{1}{r_0}>\frac{r_0-2m}{2mr_0}$. For this case,
we will consider the solution in the range $\frac{1}{r_0}\leqslant
u<\frac{1}{2m}$ and make the following substitution
\begin{equation}\label{82}u=\frac{r_0-2m}{2mr_0}+\frac{4m-r_0}{2mr_0}\sec^2\frac{\xi}{2},\end{equation}
where $u=\frac{1}{r_0}$ when $\xi=0$. By this substitution, the
equation (\ref{76}) becomes
\begin{equation}\label{83}\left(\frac{d\xi}{d\alpha}\right)^2=\frac{2m}{r_0}\left(1-k^2\sin^2\frac{\xi}{2}\right),\quad
k^2=\frac{r_0-2m}{2m}\;\;(0<k^2<1).\end{equation} Then
$\alpha$ can be expressed by the Jacobian elliptic integral as
\begin{equation}\label{85}\alpha=2\left(\frac{r_0}{2m}\right)^{\frac12}F\left(\frac{\xi}{2},k\right)\pm\frac{\sqrt{2m(r_0-2m)}}{r_0^2}\vartheta.\end{equation}

Thus, we have finished the discussions for this example.
$\qquad\qquad\qquad\blacksquare$

\section{Concluding remarks}
This paper proposes an effective  method to study the light-like
extremal surfaces in curved spacetimes, namely, one can carry out
the analysis  for light-like extremal surfaces by considering a
Cauchy problem  associated with a constraint (\ref{222}) on the
initial data. By developing a diffeomorphic mapping, this kind of
Cauchy problem for light-like extremal surfaces can be transformed
into an initial value problem for a kind of geodesic equations.
Interestingly, by this method we succeed in obtaining many
light-like extremal surfaces in Schwarzschild spacetime. Some of
them are similar to the known results in null string theory and
others are very new to our best knowledge.

As in \cite{kzz}, we believe that the equations (\ref{2}) can be
used to investigate the time-like extremal surfaces and space-like
surfaces as well as light-like surfaces. The results in this paper
verify this statement in some sense. At the same time, the Burgers
equation (\ref{10}) plays an important role in our analysis and can
be viewed as the limiting case by comparing with the key equations
(3.18) in \cite{kzz}. Finally, it will be interesting to apply the
method in this paper to consider the light-like surfaces in Kerr
spacetime or other important spacetimes.

%\vskip 10mm\noindent{\Large {\bf Acknowledgements.}} The authors
%would like to thank referees for valuable suggestions.


\begin{thebibliography}{100}
\bibitem{ac} Aurilia A. and  Christodoulou D., J. Math. Phys. {\bf 20(7)}(1979)
1446.

\bibitem{ba} Barbashov B.M., Nesterenko V.V. and Chervyakov A.M., Comm. Math. Phys. {\bf 84} (1982)
471.

\bibitem{bh} Bordemann M. and Hoppe J., Phys. Lett. {\bf B 325}(1994)
359.

\bibitem{c} Calabi E. 1970 {\it Examples of Bernstein
problems for some nonlinear equations}, in Global Analysis (Proc.
Sympos. Pure Math., Vol. XV, Berkeley, Calif., 1968), Amer. Math.
Soc., Providence, R.I., 223-230.

\bibitem{chan}Chandrasekhar S.,  {\it The mathematical theory of black
holes}, Oxford university Press Inc., New York, 1992.

\bibitem{cy} Cheng S.Y. and Yau S.T.,  Ann. of
Math. {\bf 104} (1976) 407.

\bibitem{dl}Dabrowski M.P. and Larsen A.L.,  Physical Riview D {\bf55} (1997)
6409.2

\bibitem{gordon}Gordon W. B.,  Amer. Math. Monthly {\bf79} (1972)
755.

\bibitem{go}Gorkavyy V.,  Diff. Geom. Appl. {\bf 26} (2008) 133.


\bibitem{g1} Gu C.H., {\it On the motion of a string in a curved
space-time}, Proc. of 1982 Grossmann Symposium (1982) 139-142.

\bibitem{g3} Gu C.H., {\it Extremal surfaces of mixed type in Minkowski space
$\mathbb{R}^{n+1}$.} in Variational methods (Paris, 1988), 283-296,
Progr. Nonlinear Differential Equations Appl. {\bf 4},
Birkh$\ddot{a}$user Boston, Boston, MA, 1990.

\bibitem{g4} Gu C.H.,  Chin. Ann. Math. {\bf 15B} (1994)
385.


\bibitem{hek1}He C.L. and Kong D.X., arxiv: 1002.1357v2 (2010).

\bibitem{hek}He C.L. and Kong D.X.,  arxiv: 1007.4232v1 (2010).

\bibitem{he}He Y.J.,  J. Part.  Diff. Eq. {\bf
23(2)} (2010) 158.

\bibitem{h} Hoppe J.,   Phys.
Lett. B {\bf 329} (1994) 10.

\bibitem{hk}Huang S.J. and Kong D.X., J.
Math. Phys. {\bf 48} (2007) 083510.

\bibitem{kl}Karlhede A. and Lindstr\"{o}m U,  Class. Quantum Grav. {\bf3}
(1986) L73.

\bibitem{ksz} Kong D.X., Sun Q.Y. and Zhou Y., J. Math. Phys. {\bf 47} (2006)
013503.

\bibitem{kz} Kong D.X. and Zhang Q.,  Physica D:
Nonlinear Phenomena {\bf 238} (2009) 902.

\bibitem{kzz} Kong D.X., Zhang Q. and Zhou Q., Comm.
Math. Phys. {\bf 269} (2007) 153.

\bibitem{li} Li T.T.  {\it Global classical solutions for quasiliear
hyperbolic systems}. Research in Applied Mathematics {\bf 32},
J.Wiley/Masson, 1994.

\bibitem{lin} Lindblad H.,
Proc. Amer. Math. Soc. {\bf 132} (2004) 1095.

\bibitem{m} Milnor T.,  Michigan Math. J. {\bf 37} (1990) 163.

\bibitem{schild}Schild A., Physical Review
D {\bf 16} (1977) 1722.


\end{thebibliography}
 \end{document}